\documentclass[runningheads]{llncs}

\usepackage[breaklinks, colorlinks, linkcolor=blue, citecolor=blue, urlcolor=black]{hyperref}

\usepackage{graphicx}
\usepackage{tabularx}
\usepackage{tabularcalc}
\usepackage{booktabs}
\usepackage{multirow}
\usepackage{algorithm}
\usepackage{algpseudocode}
\usepackage{tikz}
\usepackage{forest}
\usepackage{pgfplots}
\usepackage{xltabular}
\usetikzlibrary{arrows.meta}
\usepackage{amsmath}
\usepackage{amssymb}
\usepackage{xspace}
\usepackage{xcolor}
\usepackage{cleveref}

\definecolor{pbblue}{RGB}{68,119,170}
\definecolor{pbcyan}{RGB}{102,204,238}
\definecolor{pbgreen}{RGB}{34,136,51}
\definecolor{pbred}{RGB}{238,102,119}

\newcommand{\solver}[1]{\textsc{#1}\xspace}
\newcommand{\scip}{\solver{SCIP}}

\newcommand{\gurobi}{\solver{Gurobi}}
\newcommand{\soplex}{\solver{SoPlex}}

\newcommand{\papilo}{\solver{PaPILO}}
\newcommand{\papiloversion}{2.0}
\newcommand{\papilov}{\solver{PaPILO}~\papiloversion\xspace}
\newcommand{\highs}{\solver{HiGHS}}
\newcommand{\pdlp}{\solver{PDLP}}
\newcommand{\ortools}{\solver{OR-Tools}}
\newcommand{\miplibbenchmark}{MIPLIB\,2017\,Benchmark\,Set}
\newcommand{\miplibcollection}{MIPLIB\,2017\,Collection}
\newcommand{\coefftightening}{\texttt{Coeff\-Tightening}}
\newcommand{\colsingleton}{\texttt{Col\-Single\-ton}}
\newcommand{\dualfix}{\texttt{Dual\-Fix}}
\newcommand{\dualinfer}{\texttt{Dual\-Infer}}
\newcommand{\doubletoneq}{\texttt{Double\-To\-N\-Eq}}
\newcommand{\domcol}{\texttt{Dom\-Col}}
\newcommand{\implint}{\texttt{Impl\-Int}}
\newcommand{\fix}{\texttt{Fix\-Continuous}}
\newcommand{\probing}{\texttt{Probing}}
\newcommand{\parallelcols}{\texttt{Parallel\-Cols}}
\newcommand{\parallelrows}{\texttt{Parallel\-Rows}}
\newcommand{\propagation}{\texttt{Propagation}}
\newcommand{\simpleprobing}{\texttt{Simple\-Probing}}
\newcommand{\simplifyineq}{\texttt{Simplify\-Ineq}}
\newcommand{\sparsify}{\texttt{Sparsify}}
\newcommand{\stuffing}{\texttt{Stuffing}}
\newcommand{\substitution}{\texttt{Sub\-stitution}}

\newcommand{\immediately}{\texttt{1-imme\-diately}}
\newcommand{\delayed}{\texttt{1-delayed}}

\newcommand{\bracketzerozero}{$\geq 0.01$\,sec}
\newcommand{\bracketzero}{$\geq 0.1$\,sec}
\newcommand{\bracketone}{$\geq 1$\,sec}
\newcommand{\bracketten}{$\geq 10$\,sec}
\newcommand{\brackethundred}{$\geq 100$\,sec}
\newcommand{\myorcidlink}[1]{\,\href{https://orcid.org/#1}{\raisebox{-0.45ex}{\includegraphics[width=1.8ex]{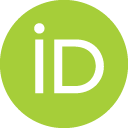}}}}

\usepackage{natbib}
\bibpunct[, ]{(}{)}{,}{a}{}{,}%
\def\bibfont{\small}%

\usepackage{todonotes,xspace}
\begin{document}

\title{\papilo: A Parallel Presolving Library for Integer and Linear Optimization with Multiprecision Support}
\author{Ambros Gleixner\inst{1,2}, Alexander Hoen\inst{1}, Leona Gottwald\inst{1}}
\institute{%
  Zuse Institute Berlin, Takustr.~7, 14195 Berlin, Germany\\
  \email{\{gleixner, hoen, gottwald\}@zib.de}
  \and
  HTW Berlin, Treskowallee 8, 10313 Berlin, Germany\\
}

\authorrunning{Gleixner, Hoen, Gottwald} 
\titlerunning{\papilo: Parallel Presolve for Integer and Linear Optimization} 
\maketitle

\abstract{Presolving has become an essential component of modern MIP solvers
	both in terms of computational performance and numerical robustness. In this
	paper we present \papilo \citep{Repository}, a new C++ header-only library that provides a large
	set of presolving routines for MIP and LP problems from the literature. The
	creation of \papilo was motivated by the current lack of (a)
	solver-independent implementations that (b) exploit parallel hardware, and (c)
	support multiprecision arithmetic.  Traditionally, presolving is designed to
	be fast.  Whenever necessary, its low computational overhead is usually
	achieved by strict working limits.  \papilo's parallelization framework aims
	at reducing the computational overhead also when presolving is executed more
	aggressively or is applied to large-scale problems.  To rule out conflicts
	between parallel presolve reductions, \papilo uses a transaction-based
	design. This helps to avoid both the memory-intensive allocation of multiple
	copies of the problem and special synchronization between
	presolvers. Additionally, the use of Intel's TBB library aids \papilo to
	efficiently exploit recursive parallelism within expensive presolving routines
	such as probing, dominated columns, or constraint sparsification. We provide
	an overview of \papilo's capabilities and insights into important design
	choices.}

\subsubsection*{Acknowledgements}
~\\
The work for this article has been partly conducted within the Research Campus MODAL funded by the German Federal Ministry  of Education and Research (BMBF grant number 05M14ZAM).

\tikzset{
	tree node/.style = {align=center, inner sep=0pt, font = \scriptsize},
	S/.style = {draw, circle, minimum size = 8mm, top color=white, bottom color=white},
	tree node label/.style={font=\scriptsize},
}
\forestset{
	declare toks={left branch prefix}{},
	declare toks={right branch prefix}{},
	declare toks={left branch suffix}{},
	declare toks={right branch suffix}{},
	tree node left label/.style={
		label=170:#1,
	},
	tree node right label/.style={
		label=10:#1,
	},
	maths branch labels/.style={
		branch label/.style={
			if n=1{
				edge label={node [left, midway] {$\forestoption{left branch prefix}##1\forestoption{left branch suffix}$}},
			}{
				edge label={node [right, midway] {$\forestoption{right branch prefix}##1\forestoption{right branch suffix}$}},
			}
		},
	},
	text branch labels/.style={
		branch label/.style={
			if n=1{
				edge label={node [left, midway] {\foresteoption{left branch prefix}##1\forestoption{left branch suffix}}},
			}{
				edge label={node [right, midway] {\forestoption{right branch prefix}##1\forestoption{right branch suffix}}},
			}
		},
	},
	text branch labels,
	set branch labels/.style n args=4{%
		left branch prefix={#1},
		left branch suffix={#2},
		right branch prefix={#3},
		right branch suffix={#4},
	},
	set maths branch labels/.style n args=4{
		maths branch labels,
		set branch labels={#1}{#2}{#3}{#4},
	},
	set text branch labels/.style n args=4{
		text branch labels,
		set branch labels={#1}{#2}{#3}{#4},
	},
	branch and bound/.style={
		/tikz/every label/.append style=tree node label,
		maths branch labels,
		for tree={
			tree node,
			S,
			math content,
			s sep'+=20mm,
			l sep'+=5mm,
			thick,
			edge+={thick},
		},
		before typesetting nodes={
			for tree={
				split option={content}{:}{tree node left label,content,tree node right label,branch label},
			},
		},
		where n children=0{
			tikz+={
				\draw [thick]  ([yshift=-10pt, xshift=-2.5pt].south west) -- ([yshift=-10pt, xshift=2.5pt].south east);
			}
		}{},
	},
}
\newcommand{\bracket}[2]{[#1,#2]}

\section{Introduction}
\label{sec::intro}

In this paper, we present a new open-source software package that provides
efficient presolving routines for mixed integer and linear optimization problems
in a framework fit for use with parallel hardware and multiprecision arithmetic.
It operates on the general class of \emph{mixed integer programs} (MIPs) given
in the form
\begin{align*}
  \min\; & c^{T}x\\
  & L \leq Ax \leq U\\
  & \ell \leq x \leq u\\
  & x_i \in \mathbb{Z}\;\; \text{for all } i \in \mathcal{I},
\end{align*}
where $A\in \mathbb{R}^{m \times n}$, $c\in\mathbb{R}^{n}, b \in \mathbb{R}^m $, and $\mathcal{I} \subseteq \{1,\dots,n\}$ is the index set of integer decision variables.
%

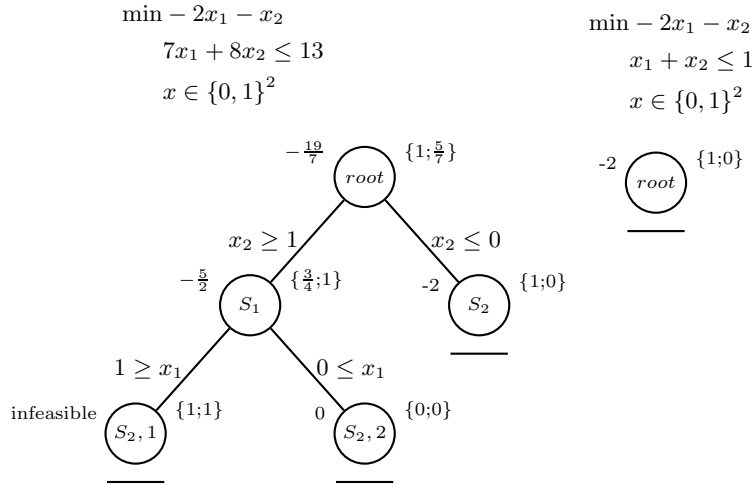
\begin{figure}[ht]
		\begin{minipage}{.48\linewidth}
			\centering
			\begin{align*}
			\min &-2x_1 -x_2\\
			& 7x_1 + 8x_2 \leq 13\\
			& x\in \{0,1\}^2 
			\end{align*}
			\begin{forest}
				branch and bound,
				where level=1{
					set branch labels={x_2\geq}{}{x_2\leq}{},
				}{
					if level=2{
						set branch labels={}{\geq x_1}{}{\leq x_1},
					}{},
				}
				[$-\frac{19}{7}$:root:\{1;$\frac{5}{7}$\}
				[$-\frac{5}{2}$:S_1:\{$\frac{3}{4}$;1\}:1
				[infeasible:{S_2,1}:\{1;1\}:1]
				[0:{S_2,2}:\{0;0\}:0]
				]
				[-2:S_2:\{1;0\}:0]
				]
			\end{forest}
		\end{minipage}
		\begin{minipage}[b]{.48\linewidth}
			\centering
			\begin{align*}
			\min &-2x_1 -x_2\\
			& x_1 + x_2 \leq 1\\
			& x\in \{0,1\}^2 
			\end{align*}
			\begin{forest}
				branch and bound,
				where level=1{
					set branch labels={x_2\geq}{}{x_2\leq}{},
				}{
					if level=2{
						set branch labels={}{\geq y}{}{\leq y},
					}{},
				}
				[-2:root:\{1;0\}]
			\end{forest}
		\end{minipage}
		\caption[Example: Comparison not presolved vs. presolved]{Comparison of the branch and bound tree of not presolved (left) vs. presolved (right).}
		\label{fig::comparison_presolve}
	\end{figure}

	
	Presolving impacts the performance of state-of-the-art MIP and LP solvers in different ways:
	By reducing the model size, it aims to accelerate the node processing time.
	Further, it often reduces the number of potential nodes in the tree to gain a speed-up for the subsequent solving process and tries to provide numerical robustness.
	On the one hand, presolving removes redundant or empty rows or columns in the model, potentially caused by imperfections in the modeling process.
	On the other hand, it uses more sophisticated techniques, for example, for tightening the linear programming (LP) relaxation or fixing variables.
	\Cref{fig::comparison_presolve} illustrates the effect of a presolving technique called coefficient tightening on a small example.
	While in the original problem, a standard branch-and-bound tree explores 5 nodes, the presolved problem requires only one solve of the LP relaxation to yield an optimal solution.
	
	%

	In one of the first and important contributions to presolving, \cite{first_presolving} introduced basic techniques like removing redundant rows, fixing variables, tightening bounds, and more.
	Over the years, many new techniques followed
	\citep{
		experiments, 
     	logical_reductions, 
     	CrowderEllisPadberg,
	  	Hoffman1991, 
	  	preprocessing1994, 
	  	presolving1995,
	  	GamrathKochMartinetal,
  	  	GemanderChenWeningeretal}. 
	In the performance analysis of \cite{12_years_progress} in CPLEX, disabling presolving results in an 11~times higher average solving time on instances that took at least 10~seconds to solve.
	In one of the most comprehensive publications to date, \cite{presolving_achterberg} describe a multitude of presolving techniques and measure their impact on the subsequent solving process within the MIP solver Gurobi. 
	In their overall performance analysis, disabling presolving worsened the solving time by the factor of 9 on instances with at least 10 seconds runtime.
   It becomes clear that presolving is often the crucial factor for whether an instance can be solved to optimality or not \citep{12_years_progress}.

	This paper presents the new software library \emph{\papilo: \textbf{Pa}rallel \textbf{P}resolving for \textbf{I}nteger and \textbf{L}inear \textbf{O}ptimization} \citep{Repository}.
    \papilo provides (a) a standalone, solver-independent implementation of
    presolving routines with (b) a natively parallel software design and (c) support
    for multi-precision computation up to exact rational arithmetic.
	To the best of our knowledge, \papilo is the first project that offers these functionalities.
	

	
   \papilo has already been used in research projects, e.g., to test the new first-order solver PDLP \citep{pdlp}.
    Since Version~7.0, the MINLP solver \scip\footnote{\url{https://www.scipopt.org/}} \citep{SCIP} integrates \papilo as one of its default presolvers. 
	Disabling the presolver plugin for \papilo{} results in a 6\% slowdown, even when called only in single-threaded mode.
    For its multiprecision features, \papilo has been particularly impactful by providing exact rational presolving for a roundoff-error-free version of SCIP \citep{EiflerGleixner2022_Acomputational}.


	%
	The paper is organized as follows.
	%
	%
	In \Cref{sec::design} we describe the algorithmic structure of \papilo, in particular its transaction-based design that is used to avoid conflicts when calling different presolving techniques in parallel.
	\Cref{sec::computational_results} presents experimental results to analyze \papilo's sequential and parallel performance.
	\Cref{sec::how-to-use} provides a brief overview how \papilo can be used as a software and \Cref{sec::conclusion} gives concluding remarks.

\section{Algorithms and software design}
\label{sec::design}

	The overall goal of presolving is to improve the performance and numerical stability of MIP solving by removing modeling artifacts, reducing and tightening the formulation, and improving its numerical properties. 
	A variety of different techniques exist, some of which may overlap or dominate each other, but many interact in a complementary fashion and rather enable each other.  
	Hence, presolving is typically organized in an iterative fashion to assure that every presolver is called at different phases of the presolving process.

	In \Cref{subsec::round_based_presolving}, we outline the overall round-based structure of \papilo's presolving loop and summarize the presolving techniques that are currently available.
In Sections~\ref{subsec::tsx_design} to \ref{subsec::internal_parallelism}, we describe how \papilo implements parallelism across and within presolvers and eliminates conflicting reductions by its transaction-based design.
Finally, \Cref{subsec::multiprecision} contains the different levels of precisions at which \papilo can perform its computations.

   \subsection{Round-based presolving in \papilo}
	\label{subsec::round_based_presolving}
	Like in most MIP solvers, \papilo organizes its presolving in rounds. Presolvers are categorized by the worst-case complexity of their implementation, per round:
	\begin{itemize}
		\item[Fast:] $\mathcal{O}(n \log n)$ with $n$ the number of non-zeros of $A$ changed since the last call.
		\item[Medium:] $\mathcal{O}(N \log N)$ with $N$ the number of non-zeros of $A$.
		\item[Exhaustive:] $\mathcal{O}(N^2)$ with $N$ the number of non-zeros of $A$.
	\end{itemize}
    A subset of the exhaustive presolvers can be marked as delayed.\footnote{Currently this is only the \sparsify{} presolver.} These are only enabled after the first time that the other exhaustive presolvers are not successful.
   
	\papilo stores and updates the maximum and minimum activity \citep{presolving_achterberg} for every constraint and therefore is able to reduce the input size of the complexity of presolvers that are primarily based on these activities. 
	\Cref{tab::implemented_presolvers} lists all presolvers currently implemented in \papilov with their complexity class.

	The presolving loop starts by executing all fast presolvers. If enough reductions are applied, it starts a new round of fast presolving regardless of the current round. Otherwise, it continues with a round at the next complexity level. The first time the non-delayed exhaustive presolvers do not find enough applicable reductions, the delayed presolvers are enabled, and \papilo proceeds with a round of fast presolvers. 
	At the beginning and after each presolving round a trivial presolving step is performed to get rid of empty columns or singleton rows and further small redundancies.

\begin{figure}[ht]
	\centering
	\scalebox{0.8}{
		\begin{tikzpicture}[node distance=2cm]
		\node (node0) [] {start};
		\node (fast) [right of=node0, draw=pbcyan] {fast};
		\node (dec1) [diamond, right of=fast, fill=pbcyan] {};
		\node (medium) [right of=dec1,  draw=pbblue] {medium};
		\node (dec2) [diamond, right of=medium, fill=pbblue] {};
		\node (exhaustive) [right of=dec2,  draw=pbgreen] {exhaustive};
		\node (dec3) [diamond, right of=exhaustive, fill=pbgreen] {};
		\node (dec4) [diamond, right of=dec3, fill=pbgreen] {};
		\node (nf) [below of=fast, fill=black, inner sep=0pt,minimum size=1pt] {};
		\node (nm) [below of=dec1,  inner sep=0pt,minimum size=1pt] {};
		\node (ne) [below of=dec2,  inner sep=0pt,minimum size=1pt] {};
		\node (nl) [below of=dec3,fill=black,  inner sep=0pt,minimum size=1pt] {};
		\node (n4) [below of=dec4,fill=black,  inner sep=0pt,minimum size=1pt] {};
		\node (end) [right of = dec4] {end};
		\node (descr1) [above of=dec1,rotate= 90, xshift=1cm] {enough reductions applied?};
		\node (descr2) [above of=dec2,rotate= 90, xshift=1cm] {enough reductions applied?};
		\node (descr3) [above of=dec3,rotate= 90, xshift=1cm] {enough reductions applied?}; 
		\node (descr4) [above of=dec4,rotate= 90, xshift=1cm] {delayed presolver disabled?};
		\node (descr4) [below of=n4, yshift = 1.5cm, xshift=1cm] {enable delayed presolvers};
		\draw[-{Stealth[length=3mm, width=2mm]}] (node0) -- (fast);
		\draw[-{Stealth[length=3mm, width=2mm]}] (fast) --  (dec1);
		\draw[-{Stealth[length=3mm, width=2mm]}] (dec1) -- node[above] {no} ++(medium);
		\draw[-{Stealth[length=3mm, width=2mm]}] (medium) -- (dec2);
		\draw[-{Stealth[length=3mm, width=2mm]}] (dec2) -- node[above] {no} ++(exhaustive);
		\draw[-{Stealth[length=3mm, width=2mm]}] (exhaustive) -- (dec3);
		\draw[-{Stealth[length=3mm, width=2mm]}] (dec3) -- node[above] {no} ++(dec4);
		\draw[-{Stealth[length=3mm, width=2mm]}] (dec4) -- node[above] {no} ++(end);
		\draw[-] (dec1.south) -- node[left] {yes} ++ (nm);
		\draw[-] (dec2.south) -- node[left] {yes} ++ (ne);
		\draw[-] (dec3.south) -- node[left] {yes} ++ (nl);
		\draw[-] (dec4.south) -- node[left] {yes} ++ (n4);
		\draw[-] (n4) -- (nf);
		\draw[-{Stealth[length=3mm, width=2mm]}] (nf) -- (fast);
		\end{tikzpicture}	
	}\caption{Workflow of PaPILO.}
	\label{fig::papilo_workflow}
\end{figure}
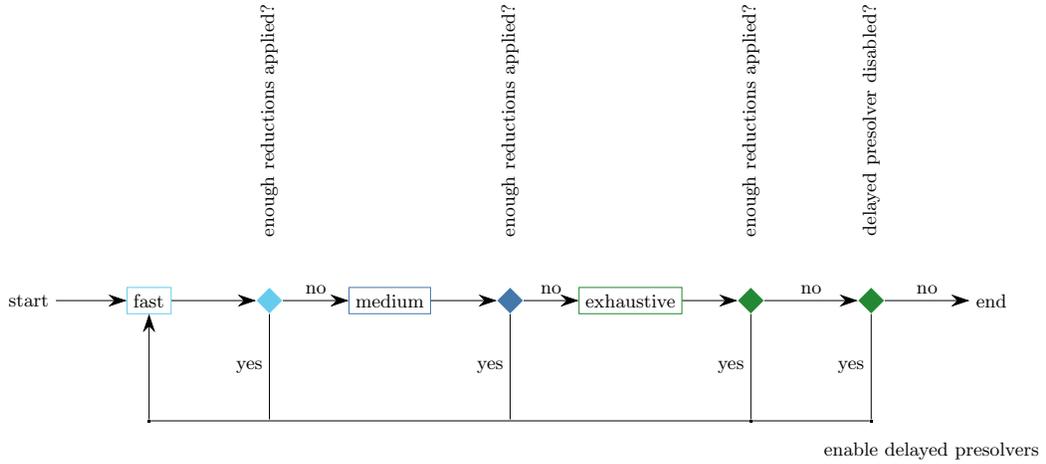

	The presolving loop terminates once the exhaustive presolvers including the exhaustive delayed presolvers do not find enough applicable reductions.\footnote{To determine, whether not enough reductions are applied the criteria a) 0.1 $\cdot$ bound changes + deleted columns $\leq\; a \cdot$ number of columns, b) side changes + deleted rows $\leq \; a \cdot $ number of rows and c) coefficient changes $\leq \; a \cdot $ number of non zeros are fulfilled. $a$ is the abort factor and can be configured via parameter settings (see \texttt{presolve.abortfac}) and manages the aggressiveness of presolving in \papilo.}
	The overall workflow is illustrated in \Cref{fig::papilo_workflow}.

	\begin{xltabular}{\linewidth}{llX}
		\hline
		Presolver & Complexity & Description \\ [0.5ex] 
		\hline\hline
		\endhead
		\coefftightening{} & fast & tightens constraint coefficients in order to strengthen the LP relaxation \citep[3.3]{presolving_achterberg}\\
		\hline 
		\propagation{} & fast & tightens variable bounds \citep[3.2]{presolving_achterberg} \\
		\hline
		\colsingleton{} & fast & removes variables appearing in only one constraint\\
		\hline
		\dualfix{} {} & medium & counts up- and down-locks on variables and exploits this information \citep[4.4]{presolving_achterberg}\\
		\hline
		\fix{} & medium & fixes continuous variables whose bounds are very close\\
		\hline
		\parallelcols{} & medium & detects parallel columns and merges them into a new variable \citep{GemanderChenWeningeretal}\\
		\hline
		\parallelrows{}{} & medium & detects parallel rows in the matrix and merges them into one \citep{GemanderChenWeningeretal}\\
		\hline
		\simpleprobing{} & medium & detects implications for binary variables on equations with special requirements \citep[3.6]{presolving_achterberg}\\
		\hline
		\doubletoneq{} & medium & searches for equations with exactly two variables and substitutes one of these variables\\
		\hline
		\simplifyineq{} & medium & deletes variables in a constraint that will never contribute to the outcome of the constraint \citep{Weninger2016}\\
		\hline
		\stuffing{} & medium & removes variables appearing in only one constraint if \colsingleton{} fails to remove them \citep{GamrathKochMartinetal}\\
		\hline
		\domcol{}  & exhaustive & detects dominated variables \citep[6.4]{presolving_achterberg}\\
		\hline
		\dualinfer{} & exhaustive &exploits complementary slackness and derives bound changes or fixes on variables \citep[7.5]{presolving_achterberg}\\
		\hline
		\implint{}  & exhaustive & detects continuous variables that are implied integers \citep[7.6]{presolving_achterberg}\\
		\hline
		\probing{} & exhaustive\footnote{The iterated propagation in \probing{} can in theory run for infinite time and hence, the complexity of \probing{} is greater than $\mathcal{O}(N^2)$. For practical reasons, \probing{} is executed together with the exhaustive presolvers.} & finds implications by probing on binary variables \citep[7.2]{presolving_achterberg}\\
		\hline
		\sparsify{} & exhaustive & adds up constraint to reduce the amount of nonzeros \citep[5.3]{presolving_achterberg} \\
		\hline
		\substitution{} & exhaustive & substitutes equations\\
		\hline
		\caption{List of presolvers implemented in \papilo.}
		\label{tab::implemented_presolvers}
	\end{xltabular}

	\subsection{Transaction-based conflict management}
	\label{subsec::tsx_design}

   The efficient parallelization of presolving requires the implementation of some kind of conflict management to avoid consistency problems.
   One reason is that many reduction techniques combine a row and column view of the constraint matrix~$A$, and decomposition of $A$ into completely independent subproblems is rarely possible.
	A na\"ive approach would be the use of thread synchronization. 
	In this case, the data is considered as a shared resource and every thread locks its required data.
	But in general, presolving techniques consist of very small ``atomic'' tasks. 
	Therefore, such a strategy comes at a high cost regarding synchronization and locking time and is not promising.
	
	Hence, \papilo follows a different approach and implements a transaction-based design. Instead of applying changes immediately, the results of presolvers, i.e., their reductions, are stored together with the data that their validity depends on. After all tasks are terminated, reductions are applied to the optimization problem, if they are still consistent with the modified problem.
	To summarize, the transaction-based design works as follows:
	\begin{itemize}
		\item Each presolver has read-only access to the data.
		\item After being called, each presolver returns its reductions as a set of transactions including criteria required to perform the reductions. 
		\item The \papilo core processes the transactions one by one after all presolvers (of the current round) are finished. It iterates through every transaction found by the presolvers in a deterministic order:
		\begin{itemize}
			\item First the core checks whether the criteria for the validity of the transaction still hold for the modified problem data.
			\item If yes, the reductions are applied and the changes in the problem data are recorded.
			\item If not, the transaction and its reductions are discarded.
		\end{itemize}
	\end{itemize}

	\Cref{fig::example} gives an example for such a transaction. Here, the \colsingleton{} presolver detects that the only nonzero entry of column 210 in the coefficient matrix is in row 133. In this case, the variable can be substituted in the objective and its value is erased from the constraint, together with modifying the bounds of the constraint. This reduction can only be applied, if the bounds of the column 210 and the row 133 were not modified since calculating the reduction.
	\begin{figure}[ht]
		\begin{verbatim}
			% information to validate the reduction
			row  -9 col 210 val 0   % BOUNDS LOCKED column 210
			row 133 col  -5 val 0   % LOCKED row 133
			% the reduction 
			row -11 col 210 val 133 % SUBSTITUTE_OBJ column 210 with row 133
			row 133 col 210 val 0   % change ENTRY (133,210) to 0.0
			row 133 col  -7 val 0   % change RHS of row 133 to infinity
		\end{verbatim}
		\centering
		\caption{Commented log of a transaction of the \colsingleton{} presolver.}
		\label{fig::example}
	\end{figure}

	\subsection{Avoiding conflicts}
	\label{subsec::avoiding_conflicts}
	Discarded transactions due to conflicts with previous reductions can make recalculations necessary and may therefore affect performance negatively. 
	Hence, it is important to reduce the conflicts between or even inside presolvers as far as possible. 
	With two simple examples we will show how \papilo avoids conflicts.

	Let's consider the following constraints with $x,y,z \in \{0,1\}$:
	\begin{align}
	y + z = 1\label{eq::conflictrow1}\\
	x + 3 y + 3z \leq 4\label{eq::conflictrow2}
	\end{align}
	Two reductions can be found: \simplifyineq{} deletes $x$ from \eqref{eq::conflictrow2}, and \eqref{eq::conflictrow1} can be substituted into \eqref{eq::conflictrow2}. 
	If the Substitution is performed first, \eqref{eq::conflictrow2} is modified, and therefore the presolver \simplifyineq{} can not be applied. 
	Applying them in the opposite order, the reductions can be applied both in one step.

	Let's consider the following constraints with $x,y \in \{0,1\}$:
	\begin{align}
		3 x + 3 y \leq 4\label{eq::parallelrow1}\\
		6 x + 6 y \geq 4\label{eq::parallelrow2}\\
		3 x + 3 y \geq 3\label{eq::parallelrow3}
	\end{align}
	Their coefficient vectors are linearly dependent.
	If we look at them pairwise, updates of the left-hand side or right-hand side might be necessary and cause a modification of the row. 
	This modification prevents \papilo from applying the second pair of parallel constraints and therefore two parallel constraints will be remaining.
	By performing this transaction in just one transaction this conflict can be avoided.

	However, note that not every conflict implies that a valid reduction is delayed. Some conflicts occur because the reductions were also found by another presolver and are already applied. 
	Or the row and column were already removed from the model and hence the reduction is negligible.

	\subsection{Parallelism across presolvers}
	\label{subsec::presolver_parallelism}

\begin{figure}[ht]
		\centering
		\begin{tikzpicture}[node distance=2cm]
		\node (core) [draw=black] {\papilo core};
		\node (p1) [right of=core, draw=black, xshift = 2cm, yshift =4cm] {Presolver 1};
		\node (p2) [below of=p1, draw=black, yshift=1cm] {Presolver 2};
		\node (p3) [below of=p2, draw=black, yshift=1cm] {Presolver 3};
		\node (plabel) [above of=p1,yshift=-1cm] {Presolvers};
		\node (nodehelp) [left of=p2, fill=black, inner sep=0pt,minimum size=1pt] {};
		
		\node (t1) [right of=core, draw=black, rotate=270, yshift = 5cm] {Transaction 5};
		\node (t2) [right of=t1, draw=black, rotate=270, yshift= -1cm] {Transaction 4};
		\node (t3) [right of=t2, draw=black, rotate=270, yshift = -1cm] {Transaction 3};
		\node (t4) [right of=t3, draw=black, rotate=270, yshift = -1cm] {Transaction 2};
		\node (t5) [right of=t4, draw=black, rotate=270, yshift = -1cm] {Transaction 1};
		\node (tlabel) [right of=t5, rotate=270, yshift = -1cm] {transactions};
		\node (tlabel) [above of=t3, yshift = 0 .5cm] {store transactions};
		\node (nodehelp2) [below of=t3, fill=black, inner sep=0pt,minimum size=1pt] {};
		\node (tlabel) [below of=nodehelp2, yshift = 1cm, align=left] {after all presolvers are finished:\\ validate transactions \textbf{sequentially}\\and apply them \textbf{conditionally}};

    	\draw[thick,dotted] ($(p1.north west)+(-0.3,0.4)$)  rectangle ($(p3.south east)+(0.3,-0.4)$);
	    \draw[thick,dotted] ($(t1.north west)+(-0.8,0.3)$)  rectangle ($(t5.south east)+(0.8,-0.3)$);
	   	\draw[thick,dotted] ($(t3.north west)+(-0.8,0.3)$)  rectangle ($(t3.south east)+(0.8,-0.3)$);
	   	\draw[-] (core)|-node[above] {start presolvers} ++(nodehelp);
	   	\draw[-{Stealth[length=3mm, width=2mm]}] (nodehelp)|-(p1);
   	   	\draw[-{Stealth[length=3mm, width=2mm]}] (nodehelp)|-(p2);
 		\draw[-{Stealth[length=3mm, width=2mm]}] (nodehelp)|-(p3);
 		\draw[-, dashed] (core)|-(nodehelp2);
 		\draw[-{Stealth[length=3mm, width=2mm]}, dashed] (nodehelp2)-|(t1.east);
 		\draw[-{Stealth[length=3mm, width=2mm]},dashed] (nodehelp2)-|(t2.east);
 		\draw[-{Stealth[length=3mm, width=2mm]},dashed] (nodehelp2)-|(t3.east);
 		\draw[-{Stealth[length=3mm, width=2mm]},dashed] (nodehelp2)-|(t4.east);
 		\draw[-{Stealth[length=3mm, width=2mm]},dashed] (nodehelp2)-|(t5.east);
	    \draw[-{Stealth[length=3mm, width=2mm]},dashed] (p1)-|(t5);
	   	\draw[-{Stealth[length=3mm, width=2mm]}] (p1)-|(t4);
	   	\draw[-{Stealth[length=3mm, width=2mm]}] (p2)-|(t3);
	   	\draw[-{Stealth[length=3mm, width=2mm]}] (p3)-|(t2);
	   	\draw[-{Stealth[length=3mm, width=2mm]}] (p3)-|(t1);
		\end{tikzpicture}
		\caption{Process of starting and applyling presolvers.}
		\label{fig_tsx_flow}
	\end{figure}
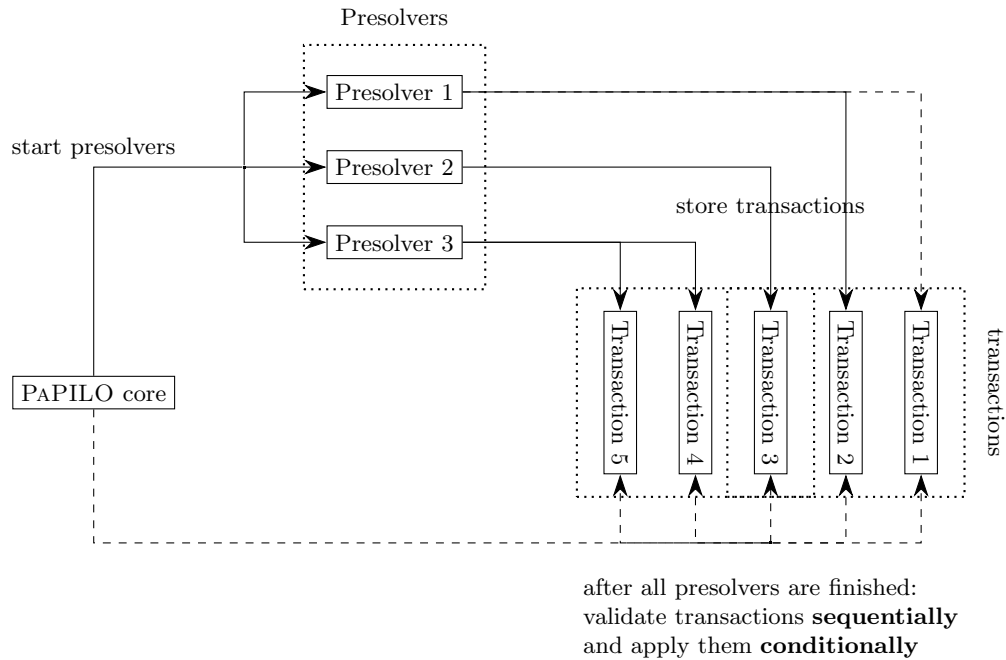

   \Cref{fig_tsx_flow} visualizes the transaction-based design described in \Cref{subsec::tsx_design}.
   Since every presolver stores its reductions in a separate data structure, 
   presolvers can be executed in parallel without any expensive synchronization.
  	\papilo relies on Intel's TBB library\footnote{\url{https://www.intel.com/content/www/us/en/developer/tools/oneapi/onetbb.html}} to execute the different tasks in parallel. 
	Since presolvers can be considered as independent tasks running on the same data, \papilo exploits task parallelism in this case.
   
   One advantage of this design is that presolvers can be implemented sequentially and will automatically participate in the parallel framework.
   Furthermore, the presolving process and result is deterministic, independently of the number of threads used.  The reason is that the reductions are stored and applied always in the same order.
   However, in the situation that only one thread is available and presolvers are necessarily executed sequentially, it may be beneficial to apply the reductions to the problem data immediately after a presolver finished and before the next one starts. This can be achieved by setting the parameter \texttt{presolve\-.apply\-\_results\-\_immediately\-\_if\-\_run\-\_sequentially} to true.


	\subsection{Parallelism within presolvers}
	\label{subsec::internal_parallelism}

   Parallelizing across presolvers only as described above can also come with potential drawbacks.
   The number of threads that can be utlized is limited by the number of presolvers per round. 
   Furthermore, the threads may not be utilized equally, since not every presolver takes the same amount of time. 

   To exploit the number of threads most efficiently, some presolving techniques in \papilo are parallelized internally.
   A particularly useful example of internal parallelization is the exhaustive presolver \probing{}. 
   We will briefly describe probing and how it is parallelized.
   Probing tries to find implications on binary variables. 
   This is done by setting a binary variable to 0 and/or 1. 
   If it turns out that the problem becomes infeasible when setting the variable to 0, the variable can be fixed to 1 in consequence; if this fixing is also infeasible, presolving can be aborted.
   Other possible implications of probing are finding new global bounds on other non-binary variables or substituting variables.

   Probing on variables are independent tasks and can therefore be parallelized.\footnote{Note that the results of probing on a single variable may still overlap, since the implication $x= 1-y$ may be found by probing on the variables $x$ and $y$. In this case, parallelization finds duplicate reductions, the latter of both will eventually be discarded.}
   Probing on a variable needs some extra storage to compute and modify the maximum and minimum activities. 
   \papilo assigns a vector to every thread, where the activities are stored. After a probing run, these vectors can be reset efficiently by reverting the changes.
   In this case \papilo exploits data parallelism because different tasks are run on the same data. 
   \papilo relies on the Intel TBB library to assign the task from the parallelism across and within the presolvers to different threads. 
   

	Further parallel presolvers in \papilo are the \domcol{}  presolver and the \sparsify{} presolver. 
	Here, the column and row hashes to detect parallel columns and parallel rows are also parallelized.
	Additionally, it can be configured via parameter settings whether the presolvers \propagation{}, \coefftightening{}, \dualfix{}, \implint{}, \simpleprobing{}, \doubletoneq{}, and \simplifyineq{} are parallelized internally. 
	In this case the iteration over the rows and variables, respectively, is parallelized.

	Certainly, not every presolver can benefit from internal parallelization. The requirement for \simpleprobing{}, e.g., are two conditions that are very fast to check. Therefore a parallelization over all constraints would result in increased runtime.

   \subsection{Multiprecision support}
   \label{subsec::multiprecision}
	Finally, \papilo provides a template for the arithmetic type used for numerical input data. 
	The template is defined as \texttt{REAL} in \papilo. 
	For the implementation of the template, \papilo relies on the Boost multi-precision package. 
	Currently, the \papilo API supports arithmetic types for double precision (default), quad precision, and rational computation with unlimited precision. 

	In the rational version all computations are performed in exact arithmetic and consequently the tolerance to compare numbers and the feasibility tolerance can be set to zero.
	We refer to the computational analysis by \cite{EiflerGleixner2022_Acomputational} for further details on the effect of exact presolving in an exact MIP framework.
   To summarize, with \papilo as presolver, more
	instances can be solved and the reduction of average solving time lies between 35.8\% and 84.1\% depending on the test set.
   The slowdown due to rational arithmetic can be partially compensated by \papilo's parallelization.
   The amount of reductions is virtually identical on MIPs with reasonable numerics, whereas for numerically difficult instances exact presolving proves to be more conservative, i.e., it yields fewer reductions.

\section{Computational results}
\label{sec::computational_results}

In this section, we analyze \papilo's performance empirically. 
In \Cref{subsec::parallel_performance} we evaluate its scalability from 1 to 32 threads. 
In \Cref{subsec::sequential_performance} we compare its sequential performance in two different modes.
In \Cref{subsec::conflict_analysis} we measure the number of conflicts appearing between individual presolvers.
The common experimental setup for our experiments is described in \Cref{subsec::experimental_setup}.

\subsection{Experimental setup}
\label{subsec::experimental_setup}

\paragraph{Testset}
We base our experiments on the extensive \miplibcollection{} \citep{miplib2017} consisting of 1,065~instances. 
Since \papilo does currently not support indicator variables or indicator constraints, instances containing such are excluded. 
Further, we exclude the large-scale instance \texttt{hawaiiv10-130}, since loading the 685,130~variables and 1,388,052~constraints of this instance exceeds the available memory.
This results in a test set of 1,022 instances. 
Every instance is presolved with two different seeds for the random number generator (except for \Cref{subsec::conflict_analysis}). 
Each seed combination is treated as a separate observation resulting in a total of 2,044~runs.

For the analysis in \Cref{subsec::conflict_analysis}, we use the detailed verbosity to log all transactions and conflicts occurring for every instance. For the \miplibcollection{} these log files become very large, therefore we restrict the experiments on the \miplibbenchmark{} resulting in an overall test set of 240~instances.

\paragraph{\papilo settings}
We vary the parameter \texttt{presolve.threads} in order to assign the number of threads to \papilo. 

If \papilo is run sequentially there is the option to apply the results after the presolver is finished instead of the end of the round. Setting the parameter \texttt{presolve\-.apply\-\_results\-\_immedia\-tely\-\_if\-\_run\-\_sequentially} to true enables this feature (\Cref{subsec::round_based_presolving}).
We refer to this setting as \immediately{}. 
If the setting is disabled we refer to this setting as \delayed.
This is also the default setting in \Cref{subsec::parallel_performance}.

To log information about transactions and conflicts in \Cref{subsec::conflict_analysis} the verbosity (\texttt{mess\-age.\-verbosity}) is set to 4 (detailed).

\paragraph{Time measurement}
The reported presolving time does not include the reading time for the instance file.
Times are always given in seconds.
For all aggregations we use the shifted geometric mean with a shift of 1.

\paragraph{Software \& Hardware}
For the experiments we use \papilo~2.0 (githash \texttt{12cd349}). 
\papilo is built with \texttt{GMP~6.1.2}\footnote{\url{https://gmplib.org/}}, \texttt{Boost~1.72}\footnote{\url{https://www.boost.org/}} and \texttt{oneTBB\-~2021.3.0}\footnote{\url{https://github.com/oneapi-src/oneTBB}}.
The experiments were carried out on identical machines with Intel(R) Xeon(R) CPU E7-8880 v4 @ 2.20\,GHz and were assigned 100,000\,MB of memory. 

\subsection{Parallel performance}
\label{subsec::parallel_performance}

In the first experiment, we want to analyze the performance of \papilo using different threads. 
The results are summarized in \Cref{tbl:papilo_performance_by_threads}.
The ``all'' row contains all instances of the test set. 
The ``$\geq n$\,sec'' row contains all instances for which at least one of the settings took at least $n$ seconds to presolve the problem; in most cases, this was the setting using just one thread.

Regardless of the number of threads used all runs return the same reductions and hence the same presolved problem.

\begin{table}[ht]
	\small
	\begin{tabular*}{\textwidth}{@{}l@{\;\;\extracolsep{\fill}}rrrrrrr@{}}
		\toprule
		& instances & 1 thread & 2 threads & 4 threads & 8 threads & 16 threads & 32 threads  \\
		\midrule
		all	& 2044 & 0.74 & 0.57 & 0.48 & 0.42 & 0.40 & 0.39 \\
		\bracketzerozero & 1753 & 0.90 & 0.69 & 0.58 & 0.51 & 0.47 & 0.46 \\
		\bracketzero & 940 & 2.22 & 1.60 & 1.28 & 1.13 & 1.01 & 0.97 \\
		\bracketone & 427 & 8.01 & 5.24 & 3.85 & 3.12 & 2.73 & 2.53 \\
		\bracketten & 134 &54.87 & 30.94&20.34 &14.79 &11.83 &10.23 \\
		\brackethundred & 49  &236.84&132.20&82.45 &55.61 &41.66 &33.93 \\
		\bottomrule
	\end{tabular*}
	\caption{Performance comparison for \papilo on \miplibcollection.}
	\label{tbl:papilo_performance_by_threads}
\end{table}

The majority of the runs, 1,617 instances, can be presolved in less than one second using even only one thread.
This results in a) an average presolving time of less than one second for all runs and b) leaves very little room for improving the overall performance by parallelization.
Note that, parallelizing (too) small chunks of work may even result in a performance decrease since the overhead of parallelization may be bigger than the benefit.
In this sense, it is a good result that the average runtime over all instances keeps decreasing up to 32 threads and almost halves the time compared to the sequential run. 

On the instances taking at least ten seconds to presolve (``\bracketten''), using 32 threads is 5.4 times faster than using one thread. For the 49 runs taking at least a hundred seconds to presolve (``\brackethundred'') the impact is even more pronounced: using 32 threads is 7.0 times faster than using one thread and 1.7 times faster than using 8 threads. In both subsets using 32 threads is 1.2 times faster than using 16 threads.

According to Amdahl's law parallelisation is limited by the percentage of sequential code.
In \papilo, this is particularly the case for modifications in the data structure like applying reductions to the matrix.
Obviously, increasing the number of threads does not affect these parts of the code. Therefore, the time for these sequential parts can be considered as fixed cost. The higher the percentage of the sequential code, the faster the speed-up factor decreases.

\begin{figure}[ht]

\begin{minipage}[t][7cm][c]{0.55\textwidth}
	\scalebox{1.0}{
\begin{tikzpicture}[scale=1,baseline]
\begin{axis}[
xlabel=threads,
xmin = 0.5, xmax = 32.5,
ymin = 0.5, ymax = 32.5,
axis lines=left,
xtick={1,2,4,8,16,32},
ticklabel style = {font=\small},
ytick={1,2,4,8,16,32},
legend cell align = {left},
legend style={nodes={scale=0.8, transform shape}}, 
legend pos = north west,
ylabel=speed-up]
\addplot[color=pbgreen, dashed, thick] coordinates {
	(1,1)
	(2,2)
	(4,4)
	(8,8)
	(16,16)
	(32,32)
};
\addplot[color=pbred, ultra thick] coordinates {
	(1,1)
	(2,1.88)
	(4,3.53)
	(8,6.2)
	(16,9.5)
	(32,12.4)
};

\addplot[domain=1:32, samples=32, color= pbcyan, dashed, thick]{1/(0.05 + 0.95/x)};
\legend{linear speed-up,  actual speed-up, speed-up for 5\% sequential code}
\end{axis}
\end{tikzpicture}
}
\end{minipage}
	\quad
\begin{minipage}[t][7cm][c]{0.3\textwidth}
	\begin{tabular}{ccc}
	threads used & time & speed-up \\
	\toprule
	1 & 102.7 & 1    \\ 
	2  & 54.6 & 1.88 \\ 
	4  & 29.1 & 3.53 \\ 
	8  & 16.5 & 6.2  \\ 
	16 & 10.8  & 9.5 \\ 
	32 & 8.3  &12.4  \\ 
	\bottomrule
\end{tabular}
\end{minipage}
	\caption{Performance of \papilo on the instance \texttt{ex10}.}
		\label{fig::probing}
\end{figure}
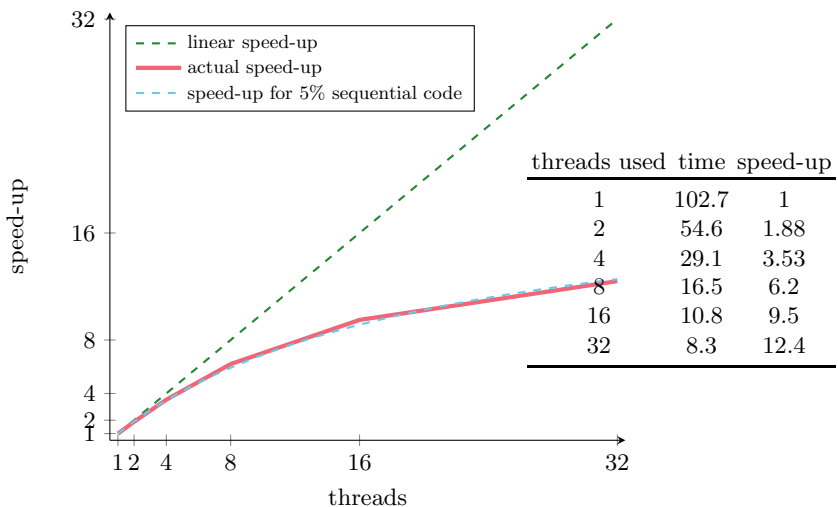

\begin{figure}[ht]
	\centering
	\begin{tikzpicture}[scale=1,baseline]
	\begin{axis}[
	xlabel=threads,
	xmin = 0.5, xmax = 32.5,
	ymin = 0.5, ymax = 10,
	axis lines=left,
	xtick={1,2,4,8,16,32},
	ticklabel style = {font=\small},
	ytick={1,2,4,8},
	legend cell align = {left},
	legend style={nodes={scale=0.8, transform shape}}, 
	legend pos = north west,
	ylabel=speed-up]
	\addplot[color=pbred, ultra thick] coordinates {
		(1,1)
		(2,1.79)
		(4,2.87)
		(8,4.26)
		(16,5.69)
		(32,6.98)
	};
	\addplot[color=pbblue, ultra thick] coordinates {
		(1,1)
		(2,1.77)
		(4,2.69)
		(8,3.71)
		(16,4.64)
		(32,5.36)
	};
	
	\addplot[domain=1:32, samples=32, color= pbgreen, dashed, thick]{1/(0.1 + 0.9/x)};
	\addplot[domain=1:32, samples=32, color= pbcyan, dashed, thick]{1/(0.2 + 0.8/x)};
	\legend{actual speed-up for set \brackethundred, actual speed-up for set \bracketten, speed-up for 10\% sequential code, speed-up for 20\% sequential code}
	\end{axis}
	\end{tikzpicture}
	\caption{Performance on sets ``\bracketten'' and `` \brackethundred'' from \Cref{tbl:papilo_performance_by_threads}.}
	\label{fig::speedup}
\end{figure}
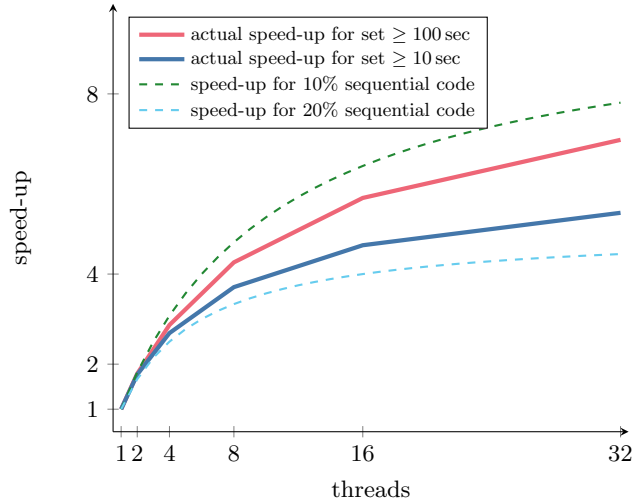

In \Cref{fig::probing} and \Cref{fig::speedup} we provide an estimation of this sequential fixed costs by comparing the actual speed-up with a linear speed-up and the speed-up for 5\%/10\%/20\% sequential code.
In \Cref{fig::probing}, we plot the presolving time of the instance \texttt{ex10} of the \miplibbenchmark{} in relation to the number of threads used.
\texttt{ex10} consists only of integer variables and the \probing{} presolver is able to reduce the size of the instance massively: from 17,860 variables and 69,608 constraints to 19 variables and 53 constraints. 
\probing{} is internally parallelized and therefore \papilo can distribute the amount of work very efficiently.

For comparison, we plot a linear speed-up and the hypothetical speed-up per thread if 95\% of the code was parallelized and the remaining five percent were run sequentially.
As can be seen, the actual speed-up on the instance \text{ex10} comes very close to the 5\% sequential code.
The overall speed-up factor achieved for 32 threads is 12.
 
\Cref{fig::speedup} shows a runtime plot for the subsets of instances that took at least 10 and 100~seconds to presolve, respectively, see \Cref{tbl:papilo_performance_by_threads}. 
Clearly, we can not expect the same results as for \texttt{ex10} where the vast majority of the time is spent in efficiently parallelized \probing{}.
Nevertheless, the speed-up for the subset with ``\bracketten'' is still slightly above the 20\% mark, indicating that 80\% of the code is executed in parallel.
For the instances taking at least $100$ seconds to presolve, the actual speed-up lies between the 10\% and 20\% sequential code mark.

\subsection{Sequential performance}
\label{subsec::sequential_performance}

One disadvantage of running presolvers in parallel is that none of the presolvers executed within the same round communicate their reductions to each other.
Therefore, \papilo has an additional sequential mode when only one thread is available (see \Cref{subsec::experimental_setup}).
Instead of applying all reductions at the end of a round, they are applied immediately after the execution of the presolver so that subsequent presolvers can already work on the modified data structure.

This change has a significant impact on the flow of \papilo and can result in a completely different sequence of rounds in \papilo.
The number of transactions found and applied per round differs, which can impact the evaluation of the next round as well as the stopping criteria at the end.

In \Cref{tbl::papilo_performance_by_1_thread}, we compare the performance of \delayed{} and \immediately{}.
Especially, we want to focus on the impact of applying the results at the end of the rounds.
Therefore, \Cref{tbl::papilo_performance_by_1_thread} contains not only the time but also the mean of the number of presolving rounds and the shifted geometric mean of the number of non-zeros of the presolved problem.

Infeasible or unbounded problems detected by \papilo are excluded from this table to be able to compare the number of non-zeros. 
Further, we exclude four instances \texttt{proteindesign*} since in \immediately{} mode a threshold is passed that triggers a very expensive probing call, while in \delayed{} presolving is stopped before probing. 
This would distort the results and will be addressed in future versions of \papilo.
All in all, this results in a test set of 990 instances.

\begin{table}[ht]
	\small
	\begin{tabular*}{\textwidth}{@{}l@{\;\;\extracolsep{\fill}}rrrrrrr@{}}
		\toprule
		& & \multicolumn{3}{c}{\immediately} & \multicolumn{3}{c}{\delayed} \\
		\cmidrule{3-5}
		\cmidrule{6-8}
		& instances & time & rounds & non-zeros & time & rounds & non-zeros \\
		\midrule
		all			&1,980 & 0.78 & 14.54 & 56,258   & 0.73 & 23.53 & 56,803   \\
		\bracketzerozero &1,728 & 0.94 & 15.98 & 96,122   & 0.87 & 26.20 & 97,190   \\
		\bracketzero & 975 & 2.13 & 22.15 & 309,525  & 1.96 & 39.09 & 309,701  \\
		\bracketone	& 457 & 7.18 & 29.15 & 692,599  & 6.62 & 59.64 & 693,300  \\
		\bracketten  & 150 &38.22 & 39.06 & 1,339,867 &40.22 &129.38 & 1,339,952 \\
		\brackethundred &  51 &184.78& 22.67 & 2,140,532 &209.91&303.67 & 2,231,813 \\
		\bottomrule
	\end{tabular*}
	\caption{Performance comparison for \papilo using \delayed{} vs \immediately{}  on \miplibcollection.}
	\label{tbl::papilo_performance_by_1_thread}
\end{table}


Enabling \immediately{} results in 38\% fewer rounds per instance compared to the normal mode (see \Cref{tbl::conflicts}) since a) no conflicts between presolvers prevent applying transactions and b) working on the updated data structure is more efficient.
The number of final non-zeros is minimally reduced by at most 3\% on ``\brackethundred''.
  
Surprisingly, despite the fewer amount of rounds \immediately{} is overall slower.
A more detailed analysis shows that the \immediately{} setting results in much fewer fast presolving (14.72 vs. 4.37) rounds but uses slightly more exhaustive rounds.
An example for this behavior is the instance \texttt{neos-1122047} where the first 547 rounds were fast rounds resulting in 7 exhaustive rounds compared to 30 exhaustive rounds for \immediately{}.
This is reducing the presolving time from 11 to 5 seconds.

On the other hand, if exhaustive presolvers play an important role during presolving and take a lot of time the \immediately{} setting is able to reduce the amount of exhaustive presolving rounds.
For example on the instances  \texttt{neos-4332801-seret}, \texttt{neos-4332810-sesia} and, \texttt{moj-mining} one round less of exhaustive presolving is necessary resulting in saving 100 seconds.
This explains why \immediately{} is slightly faster on the ``\bracketten'' and ``\brackethundred'' subsets.

\subsection{Experiments on conflicts}
\label{subsec::conflict_analysis}
Due to its transaction-based design, the parallel performance of \papilo depends to a large extent on its ability to avoid conflicts between transactions. 
Therefore, in this section, we analyze which presolvers generate conflicts and how this impacts the performance of \papilo.
We focus on the following questions:
\begin{itemize}
	\item Which presolvers interact with each other and which of these generate conflicts? (See \Cref{tbl::fast_presolver_conflicts,tbl::medium_presolver_conflicts,tbl::exhaustive_presolver_conflicts}.)
	\item How many conflicts appear and how does this impact the performance of \papilo? (See \Cref{tbl::conflicts}.)
\end{itemize}

To answer these questions we use the verbose mode of \papilo to log every transaction and its status (applied/discarded/canceled).
Please note again that we restrict these experiments to the \miplibbenchmark{} \citep{miplib2017} consisting of 240 instances, since these log files can become very large.

A necessary precondition for conflicts between two presolvers is obviously that these presolvers find transactions within the same round. 
Therefore, we analyze on how many rounds a pair of presolvers finds transactions within the same round and on how many of these rounds a conflict appears.
\Cref{tbl::fast_presolver_conflicts,tbl::medium_presolver_conflicts,tbl::exhaustive_presolver_conflicts} report these results for fast, medium and exhaustive presolvers, respectively.

\Cref{tbl::fast_presolver_conflicts,tbl::medium_presolver_conflicts,tbl::exhaustive_presolver_conflicts} are organized as follows:
The presolvers are listed in the same order in which their transactions are applied by the core. 
The first column ``calls'' indicates on how many rounds the presolver found at least one transaction. 
The remaining columns indicate on how many of their common calls the presolver in the row and the presolver in the column found both at least one pair of conflicting transactions appeared.
Note that transactions of one and the same presolver may also be in conflict.
 
 
\begin{table}[ht]
	\centering
	\begin{tabular*}{0.9\textwidth}{@{}l@{\;\;\extracolsep{\fill}}rccc@{}}
		\toprule
		\multicolumn{1}{c}{} &\multicolumn{1}{c}{} & \multicolumn{3}{c}{conflicting calls/common calls}\\
		\cmidrule{3-5}
		presolver & calls & \multicolumn{1}{c}{cs} & \multicolumn{1}{c}{co} & \multicolumn{1}{c}{cp} \\
		\midrule
		colsingleton (cs) & 505 & 45/505 & -  & - \\
		coefftightening (co) & 529 & 0/244 & 0/529 & - \\
		propagation (cp) & 1,530 & 41/384 & 0/401 & 0/1,530\\
		\bottomrule
	\end{tabular*}
	\caption{Conflicting/common calls on 1728 rounds of fast presolvers.}
	\label{tbl::fast_presolver_conflicts}
\end{table}

\begin{table}[ht]
	\centering
	\begin{tabular*}{0.9\textwidth}{@{}l@{\;\;\extracolsep{\fill}}rrccccccc@{}}
		\toprule
		\multicolumn{1}{c}{} & \multicolumn{1}{c}{} & \multicolumn{8}{c}{conflicting calls/common calls}\\
		\cmidrule{3-10}
		presolver & calls & \multicolumn{1}{c}{sp} & \multicolumn{1}{c}{pr} & \multicolumn{1}{c}{pc} & \multicolumn{1}{c}{st} & \multicolumn{1}{c}{d} & \multicolumn{1}{c}{f} & \multicolumn{1}{c}{si} & \multicolumn{1}{c}{dt} \\
		\midrule
		simpleprobing (sp) & 36 & 0/36 & -  & -  & -  & -  & -  & -  & - \\
		parallelrows (pr) & 192 & 0/12 & 0/192 & -  & -  & -  & -  & -  & - \\
		parallelcols (pc) & 41 & 0/1 & 0/17 & 0/41 & -  & -  & -  & -  & - \\
		stuffing (st) & 62 & 0/3 & 0/18 & 0/8 & 3/62 & -  & -  & -  & - \\
		dualfix (d) & 291 & 0/4 & 0/63 & 0/26 & 3/44 & 0/291 & -  & -  & - \\
		fixcontinuous (f) & 1 &  0/0  & 0/1 &  0/0  & 0/1 & 0/1 & 0/1 & -  & - \\
		simplifyineq (si) & 25 & 0/2 & 6/23 & 0/2 & 0/4 & 0/8 & 0/1 & 0/25 & - \\
		doubletoneq (dt) & 238 & 2/25 & 12/72 & 0/7 & 8/22 & 0/70 & 0/1 & 0/20 & 75/238\\
		\bottomrule
	\end{tabular*}
		\caption{Conflicting/common calls on 553 rounds of medium presolvers.}
	\label{tbl::medium_presolver_conflicts}
\end{table}

\begin{table}[ht]
	\centering
	\begin{tabular*}{0.9\textwidth}{@{}l@{\;\;\extracolsep{\fill}}rrccccc@{}}
		\toprule
		\multicolumn{1}{c}{} & \multicolumn{1}{c}{} & \multicolumn{6}{c}{conflicting calls/common calls}\\
		\cmidrule{3-8}
		presolver & calls & \multicolumn{1}{c}{ii} & \multicolumn{1}{c}{dc} & \multicolumn{1}{c}{di} & \multicolumn{1}{c}{po} & \multicolumn{1}{c}{su} & \multicolumn{1}{c}{sp} \\
		\midrule
		implint (ii) & 57 & 0/57 & -  & -  & -  & -  & - \\
		domcol (dc) & 54 & 0/6 & 20/54 & -  & -  & -  & - \\
		dualinfer (di) & 15 & 0/7 & 1/4 & 0/15 & -  & -  & - \\
		probing (po) & 862 & 0/18 & 0/12 & 0/6 & 0/862 & -  & - \\
		substitution (su) & 645 & 0/32 & 1/31 & 1/9 & 214/314 & 383/645 & - \\
		sparsify (sp) & 138 & 0/3 & 0/5 &  0/0  & 2/35 & 26/51 & 50/138\\
		\bottomrule
	\end{tabular*}
	\caption{Conflicting/common calls on 1457 rounds of exhaustive presolvers.}
	\label{tbl::exhaustive_presolver_conflicts}
\end{table}

As \Cref{tbl::fast_presolver_conflicts,tbl::medium_presolver_conflicts,tbl::exhaustive_presolver_conflicts} illustrate presolvers rarely interact conflict-wise with each other. Only 16 of the 71 pairings cause conflicts. Further, only the presolvers \colsingleton, \doubletoneq, \substitution{} and \sparsify{} have conflicts on more than 5\% percent of their rounds.

Presolvers typically require a certain problem or nonzero structure in the matrix, and thus, will not necessarily find reductions a) on the same row/columns or b) even on the same instances. 
In these cases, conflicts will not appear and \papilo, therefore, utilizes the strength of parallelization:
\papilo exploits this behavior by distributing the runtime of the (``unsuccessful'') presolver to different threads at no cost since these presolvers do not return reductions and hence can not cause conflicts.

~\\
Not only the amount of rounds with conflicts is important for performance.
The reason, type and amount of conflicts play also an important role in the performance of \papilo.
Therefore, in \Cref{tbl::conflicts} we take a closer look at the conflicts within the presolvers and analyze why they appear. 

To describe the structure of \Cref{tbl::conflicts} we introduce some notation: 
Let $N$ be the set of all instances of our test set,
let $t_i^p$ denote the number of transactions found by presolver $p$ on instance $i\in N$, and let
$c_i^{p-q}$ denote the number of transactions of presolver $p$ that could not be applied because it conflicted with a transaction of presolver $q$ on instance $i\in N$.
$r_i^{p-q}$ shall denote the number of redundant conflicts between presolver $p$ and $q$ on instance $i\in N$ (see \Cref{subsec::avoiding_conflicts}). 
Note that the technique to detect redundant conflicts is a simple heuristic and therefore $r^{p-q}$ should be considered a lower bound.
The columns in \Cref{tbl::conflicts} are organized as follows: 
\begin{itemize}
	\item The first column $p-q$ defines the pair of conflicting presolvers $p$ and $q$ for all pairs of presolvers with at least one conflict. Presolver $q$ is applied first, hence transactions of $p$ are discarded because of presolver $q$. 
	\item Column $\sum_i t_i^p$  contains the total amount of transactions found by presolver $p$ over all instances $i\in N$.
	\item Column $\sum_i c_{i}^{p-q}$ contains the total amount of conflicts between presolvers $p$ and $q$ over all instances $i \in N$.
	\item Column $\tfrac{\sum_{i: t^p_i>0} c_{i}^{p-q}/t_i^p}{|\{i| t_i^p >0 \}|}$ reports the average conflict rate on instances with at least one transaction found by presolver $p$.
	\item Column $\frac{\sum_i r_i^{p-q}}{\sum_i c_{i}^{p-q}}$ reports the percentage of redundant reductions $r_i^{p-q}$ among the conflicts between presolver $p$ and $q$.
\end{itemize}

\begin{table}[ht]
	\centering
	\begin{tabular*}{\textwidth}{@{}l@{\;\;\extracolsep{\fill}}|r|rr|r@{}}
		\toprule
		$p - q$ & $\sum_i t_i^p$ 
		 & $\sum_i c_{i}^{p-q}$ & $\tfrac{\sum_{i: t^p_i>0} c_{i}^{p-q}/t_i^p}{|\{i| t_i^p >0 \}|}$
		& $\frac{\sum_i r_i^{p-q}}{\sum_i c_{i}^{p-q}}$ \\
		\midrule
		colsingleton-colsingleton & 324,750 & 90,217 & 7.40\% & 0\% \\
		propagation-colsingleton & 23,100,091 & 7,218 & 4.03\% & 100\% \\
		stuffing-stuffing & 36,676 & 12 & 1.21\% & 8.33\% \\
		dualfix-stuffing & 233,593 & 5 & 0.02\% & 100\% \\
		simplifyineq-parallelrows & 4,734 & 288 & 14.38\% & 50\% \\
		doubletoneq-simpleprobing & 556,857 & 18 & 0.00\% & 0\% \\
		doubletoneq-parallelrows & 556,857 & 27,059 & 2.75\% & 100\% \\
		doubletoneq-stuffing & 556,857 & 177 & 1.22\% & 100\% \\
		doubletoneq-doubletoneq & 556,857 & 69,406 & 6.28\% & 0\% \\
		domcol-domcol & 128,069 & 46,789 & 14.40\% & 0\% \\
		dualinfer-domcol & 532 & 2 & 2.38\% & 100\% \\
		substitution-domcol & 247,574 & 2 & 0.00\% & 100\% \\
		substitution-dualinfer & 247,574 & 1 & 0.00\% & 0\% \\
		substitution-probing & 247,574 & 26,756 & 14.21\% & 0\% \\
		substitution-substitution & 247,574 & 81,848 & 20.58\% & 0\% \\
		sparsify-probing & 42,846 & 7 & 0.03\% & 0\% \\
		sparsify-substitution & 42,846 & 1,192 & 1.82\% & 0\% \\
		sparsify-sparsify & 42,846 & 5,766 & 7.52\% & 0\% \\

		\bottomrule
	\end{tabular*}
	\caption{Conflicts between presolvers on the \miplibbenchmark.}
	\label{tbl::conflicts}
\end{table}

Next, we describe the relationship between the presolver leading to the conflicts in \Cref{tbl::conflicts}:

\paragraph{ColSingleton and Stuffing}
	The presolvers \colsingleton{} (27.7\% of its overall transactions) and \stuffing{} (0.0\%) generate conflicts to themselves if one row contains multiple singleton columns. Currently, both presolvers iterate ``randomly'' over the cached singleton columns. Avoiding this would require a) reordering of singleton columns within the presolver or b) enriching the data structure of the singleton columns with their row. This would require further updates if rows are deleted.
	Since both presolvers are usually very fast, avoiding these conflicts would likely not improve the performance. 
		
	Further, if the \colsingleton{} presolver marks a row as redundant then the remaining singleton columns within this row become empty columns and can be deleted in the trivial presolve step and no further actions are necessary.
		
	Multiple singleton columns in a row do not occur on a regular basis in the test set, and the majority of the conflicts are generated by the instances \texttt{supportcase12} (35,704 conflicts) and \texttt{savsched1} (36,471 conflicts).

\paragraph{DoubleToNEq}
	96660 transactions of \doubletoneq{} (17.4\% of total) are discarded due to conflicts. \doubletoneq{} substitutes variables in equations with two variables. On this type of constraints,  the presolvers \simpleprobing{}, \parallelrows{}, and \stuffing{} can also find reductions leading to implications, marking the row as redundant, or relaxing the bounds. In these cases, the substitution becomes obsolete.
		
	Further, if \doubletoneq{} suggests substituting a variable that was already substituted by a previously found reduction of \doubletoneq{} this results in a conflict. This is the case if a variable appears in multiple equations with two variables. 
	In order to delete the remaining equations with two variables additional runs will be necessary. 
	These conflicts could be reduced if \doubletoneq{} would keep track of its suggested variable substitutions. But like for \colsingleton{} this will probably come at a higher cost.

\paragraph{SimplifyIneq-ParallelRows}
	288 of the 4,734 transactions of \simplifyineq{} are discarded because they are in conflict with transactions of the \parallelrows{} presolver. On parallel rows, \simplifyineq{} can find an equivalent reduction on every row and generate a separate transaction for each of these rows.
	Simultaneously, the parallel rows are detected by the \parallelrows{} presolver, marked as redundant and removed from the problem (except for one). 
	These removed rows can not be simplified anymore and the corresponding transaction is discarded. 
	Simplifications can be applied to the remaining (parallel) constraint if the bounds were not modified by \parallelrows{}. Otherwise, one additional round is necessary.

\paragraph{DomCol-DomCol}
	46,789 out of 128,069 transactions of the \domcol{} presolver are discarded due to conflicts to reductions generated by the presolver itself.
	On average, 13.3\% of these reductions are redundant. 
	If the columns $a$ and $b$ dominate column $c$ \papilo will generate two transactions - one for $a$ dominating $c$ and one for $b$ dominating $c$ - but with different requirements for the modified data structure. 
	Because of the different requirements these reductions are not considered as redundant in \Cref{tbl::conflicts}. 
	To increase the possibility that one of the two transactions can be applied both are generated.
	The majority of the conflicts is generated by the instances \texttt{co-100} (19,692 conflicts) and \texttt{nw04} (20,909 conflicts).

\paragraph{Redundant Conflicts}
	In rare cases, a presolver may find the same or weaker reductions as another presolver. For example a singleton variable can be fixed by \stuffing{} as well as by \dualfix{} (a total of 5 conflicts). Further, this happens also for \propagation-\colsingleton{} (7,218 conflicts) and \dualinfer{}-\domcol (2 conflicts).

\paragraph{Sparsify and Substitution}
	6,965 transactions of \sparsify{} and 108,607 transactions of \substitution{} are discarded.
	The vast majority because of ``internal'' conflicts. 
	The transactions of these presolvers can only be applied if the specified row is not modified. 
	This row can be modified by all presolvers applying transactions before them or even by themselves and therefore generates conflicts and additional rounds. 
	Note that \papilo cancels a ``conflict-free'' substitution if applying it increases the number of non-zeros entries.
	Consequently, discarded transactions could be possibly canceled due to this condition.

~\\
To summarize, conflicts within the same round appear rarely and in most cases, these conflicts are redundant or can be resolved with one additional round. Only the presolvers \doubletoneq{}, \substitution{} and \sparsify{}, which all perform some kind of substitution, generate a notable amount of non-redundant conflicts. 
However, the parallel speedup measured seems to outweigh some additional rounds caused by these conflicts: using 32 threads results in a speed-up of 5.4 times on instances taking at least ten seconds to presolve (see \Cref{tbl:papilo_performance_by_threads} in \Cref{subsec::parallel_performance}).

\section{How to use \papilo}
	\label{sec::how-to-use}
	As already indicated in the previous sections, \papilo can be used in three different ways: as a standalone application, as a frontend to an MIP or LP solver, and integrated into a solver via its API.
	
	When used as a standalone application, the user must solve the reduced problem with a solver of his choice. Then, \papilo can be called via the command line with the generated reduced solution to perform (dual) postsolving\footnote{The majority of the presolvers support postsolving also the dual solution, reduced costs and basis information for LP.} \citep{SCIP8} in order to obtain the (original) solution for the original problem.
	To avoid this manual postsolving step \papilo provides interfaces to MIP and LP solvers.
	By using the \texttt{solve} keyword, \papilo presolves the specified instance, then calls one of its integrated solvers via API to solve the reduced problem, and then postsolves the obtained reduced solution automatically.
	Currently, \papilo supports the MIP solvers \gurobi\footnote{\url{https://www.gurobi.com/}}, \highs\footnote{\url{https://highs.dev/}}, \scip\footnote{\url{https://scip.zib.de/}}, and the LP solver \soplex\footnote{\url{https://soplex.zib.de/}}.
	Furthermore, \ortools\footnote{\url{https://developers.google.com/optimization}} is integrated allowing to use all its available solvers such as the first-order method \pdlp \citep{pdlp}.

	Finally, \papilo also provides an interface such that it can be integrated conveniently into existing solvers.
	An example on how this is possible is available in the solver \scip \citep{SCIP}, which calls \papilo as one of its default presolving plugins.
	
	For adding new presolving methods in \papilo it is sufficient to implement the interface \texttt{PresolveMethod<REAL>} and then register the new method as a default presolver.

\section{Conclusion}
	\label{sec::conclusion}
	
	In this paper, we provided an overview of \papilo's functionality and its design pattern to allow parallel presolving. 
	
	The design goal of \papilo is to go beyond parallelizing work that is a priori independent.  Since different presolvers require different structures in the matrix in order to be applicable, they do not necessarily find reductions (a) on the same rows and columns or (b) even on the same instances.
	\papilo{}'s transaction-based design guarantees a posteriori, after attempting different presolving reductions in parallel, that no conflicting reductions are applied.
   In addition we analyzed empirically that conflicts appear rarely and are not a performance bottleneck when they occur.
	
	The parallelization is particularly useful when presolving time-consuming instances. 
	On instances taking at least 10 seconds to presolve, \papilo is able to reduce presolving time by a factor of five (\Cref{tbl:papilo_performance_by_threads}).
	But even on the entire \miplibcollection{} set, which consists mostly of instances presolved in less than one second, \papilo is able to reduce presolving time by almost 50\% in shifted geometric mean.
	 

\bibliographystyle{informs2014}
\bibliography{bibliography} 

\end{document}